\documentclass[12pt]{article}
\usepackage{latexsym}
\usepackage{amssymb,amsbsy,amsmath,amsfonts,amssymb,amscd}
\def \vp{\varphi}
\def \< {\langle}
\def \> {\rangle}
\def \t {\tilde}

\def \al{\alpha}
\def \e {\varepsilon}

\def \wh{\widehat}
\def \vp{\varphi}
\def \f {\frac}
\def \lp {{\Bigl (}}
\def \rp {{\Bigr )}}
\def \lc {\Bigl [}
\def \rc {\Bigr ]}
\def \tn {\vert \kern-.1em \vert \kern-.1em \vert}
\def\bu{\bar{u}}
\def\intrd{\int_{\RR^d}}
\def\intbr {\int_{| x |\leq R}}

\def\intbrc {\int_{| x |>R}}

\def \xox   { {x \over {\vert x \vert}} }

\def \ue    {u_\e}
\def \dr    { {\partial \over \partial r}}
\def \dtau  { {\partial \over \partial \tau}}
\def \domega  { {\partial \over \partial \omega}}
\newcommand{\fer}[1]{(\ref{#1})}
\newcommand{\beq}{\begin{equation}}
\newcommand{\eeq}{\end{equation}}

\newcommand {\commentout}[1] {{}}
\newcommand{\R}{\mathbb{R}}
\newcommand{\N}{\mathbb{N}}
\newcommand{\Z}{\mathbb{Z}}
\def \ZZ {\Z}
\def \NN {\N}
\newcommand{\Sph}{\mathbb{S}}
\newcommand {\p}   {\partial}

\newcommand {\proof} {\noindent {\bf Proof}. }
\newcommand{\beqa}{\begin{eqnarray}}
\newcommand{\bea} {\begin{array}{ll}}
\newcommand{\beqan}{\begin{eqnarray*}}
\newcommand{\eeqa}{\end{eqnarray}}
\newcommand{\eeqan}{\end{eqnarray*}}
\newcommand{\eea} {\end{array}}
\newtheorem{theorem}{Theorem}[section]
\newtheorem{lemma}[theorem]{Lemma}

\def\RR{\rm \hbox{I\kern-.2em\hbox{R}}}
\def\NN{\rm \hbox{I\kern-.2em\hbox{N}}}
\def\ZZ{\rm {{\rm Z}\kern-.28em{\rm Z}}}
\def\CC{\rm \hbox{C\kern -.5em {\raise .32ex \hbox{$\scriptscriptstyle
|$}}\kern -.22em{\raise .6ex \hbox{$\scriptscriptstyle |$}}\kern .4em}}
\def\vp{\varphi}
\def\<{\langle}
\def\>{\rangle}
\def\t{\tilde}

\def\e{\varepsilon}

\def\wh{\widehat}
\def\Chi{\raise .3ex \hbox{\large $\chi$}}
\def\vp{\varphi}
\def\[{\Bigl [}
\def\]{\Bigr ]}
\def\({\Bigl (}
\def\){\Bigr )}
\def\[{\Bigl [}
\def\]{\Bigr ]}
\def\({\Bigl (}
\def\){\Bigr )}

\def\bu{\bar{u}}
\def\intrd{\int_{\RR^d}}
\def\intbr {\int_{B(R)}}

\def\intbrc {\int_{\vert x \vert >R}}

\def \tn {\vert \kern-.1em \vert \kern-.1em \vert}

\def \trait (#1) (#2) (#3){\vrule width #1pt height #2pt depth #3pt}
\def \Box{\hfill
         \trait (0.1) (5) (0)
         \trait (5) (0.1) (0)
         \kern-5pt
         \trait (5) (5) (-4.9)
         \trait (0.1) (5) (0)
         }
\newtheorem{Tm}{Theorem}[section]
\newtheorem{Prop}[Tm]{Proposition}

\begin{document}
\centerline{\Large\bf Energy concentration and Sommerfeld condition for
Helmholtz 
%%%%%%%%%%%%%%%%%23205
}
\centerline{\Large\bf equation with variable index at infinity}
\bigskip
\bigskip
\centerline{\bf Benoit Perthame ($^*$) and Luis Vega ($^{**}$)}
\bigskip
\centerline {($^{*}$)  Ecole Normale Sup\'erieure, DMA, UMR8553}
\centerline {45, rue d'Ulm 75230 Paris, France}
\centerline {email: benoit.perthame@ens.fr}
\centerline { ($^{**}$)  Universidad del Pais Vasco, Apdo. 644}
\centerline { 48080 Bilbao, Spain}
\centerline {email: mtpvegol@lg.ehu.es}
\bigskip
\centerline{\bf Abstract}
We consider the Helmholtz equation with a variable index of
refraction $n(x)$, which is not necessarily constant at infinity but can
have an angular dependency like
$n(x)\to n_\infty(x/|x |)$ as $|x |\to \infty$. Under some appropriate
assumptions on this convergence and on $n_\infty$ we prove that the
Sommerfeld condition at infinity still holds true under the explicit form
$$  \int_{\R^d} \left| \nabla u -i n_\infty^{1/2} 
u \ \xox\, \right|^2  \ \f{dx}{|x |}<+\infty.
$$
It is a very striking and unexpected feature that the index $n_{\infty}$ appears in this
formula and not the gradient of the phase as established by Saito in  \cite {S} and broadly used numerically. This apparent contradiction is clarified by the existence of some extra
estimates on the energy decay. In particular we prove that
$$
 \int_{\R^d} \left| \nabla_\omega n_\infty(\xox)\right|^2 \, \f{ | u |^2}{|x |} \ dx <  +\infty.
$$
In fact our main contribution is to show that this can be interpreted as a concentration of the energy along the critical lines of $n_\infty$. In other words, the Sommerfeld condition hides the main physical effect arising for a variable $n$ at infinity; energy concentration on lines rather than dispersion in all directions. 
%
%
%%%%%%%%%%%%%%%%%%%%%%%%%%%%%%%%%%%%%%%%%%%%%%%%%%%
\section{ Introduction}
\setcounter{equation}{0}
%%%%%%%%%%%%%%%%%%%%%%%%%%%%%%%%%%%%%%%%%%%%%%%%%%%
%%%%%%%%%%%%%%%%%%%%%%%%%%%%%%%%%%%%%%%%%%%%%%%%%%%
%
%

We consider the Helmholtz equation with a variable index of
refraction $n(x)$, with a slow, and only radial decay to a constant
$n_{\infty}(x/|x |)$ at infinity
\beq 
i \e \ue + \Delta \ue +n(x) \ue = - f(x), \qquad \e>0.
\label{Eq:H}
\eeq
Our main interest is the so called limiting absorption principle (i.e. to
study the limit when $\e>0$  approaches to $0$ in \fer{Eq:H})
and the validity of the Sommerfeld radiation condition at infinity. One
of the main results in this paper is to prove that
\beq
  \int_{\R^d} \left| \nabla u(x)  -i n_\infty^{1/2}(\xox)  u(x) \; \xox\,\right|^2 \ \f{dx}{|x |}
<+\infty.
\label{Eq:S}
\eeq
A direct consequence of this condition is the more classical setting
\beq
{\rm liminf} \int_{|x |=r} \left| \nabla u(x) -i n^{1/2}_{\infty}(\xox) u(x) \; \xox \right|^2
\  d\sigma(x) \to 0, \qquad {\rm as}
\quad r \to \infty ,
\label{Eq:SP}
\eeq
where $d\sigma$ denotes the Lebesgue measure on the sphere.
It is a very striking and unexpected feature that the term  
$$
n_{\infty}^{1/2}(\xox) \; \frac{x}{|x|}
$$ 
appears in this formula instead of $\nabla \varphi (x)$ with $\varphi (x)$ the solution to the corresponding eikonal equation
$$
| \nabla \varphi|^2=n,
$$
as established by Saito in  \cite {S}.  This phenomenon, as well as the proof of \fer{Eq:SP},
can be explained by the existence of some  new energy estimate that we
state later on and which is, in some sense, the main result of this paper. It explains that the Sommerfeld radiation condition hides the main physical effect of a variable $n$ at infinity; energy is not dispersed in all directions but concentrated on those given by the critical points of $\nabla n(\xox)$. It would be interesting to prove that only local maxima of $n$ can carry energy. \\

These estimates use in a strong way the inequality
 obtained
in \cite{PV} for the tangential part of the gradient of the solution of \fer{Eq:H}. In order to recall
that result we need some notation. Firstly we define, for $j\in \Z$, the annulus
$C(j)$ by
$$  C(j)= \{ x\in \R{^d} \; s.t.\;  2^{j}\le | x| \le 2^{j+1}\}.
$$
Then we set
\beq\tn u \tn^{2}_{R_0} := \sup_{R>R_0} \frac{1}{R} \intbr \vert u \vert^2 \ dx,
\label{Eq:N1}
\eeq
\beq N_{R_0}(f): = \sum_{j>J} [ 2^{j+1} \int_{C(j)}
| f|^{2} \ dx ]^{1/2} + [R_0\int_{B_{R_0}}\,|f|^2 \ dx]^{1/2},
\label{Eq:N2}
\eeq
with $J$ defined by $2^J\leq R_0 <2^{J+1}$, and we drop the index $R_0$
if  $R_0=0$.
We also denote the radial and tangential derivatives by
\beq
\dr u(x):= \xox \cdot \nabla u(x), \qquad \nabla_\tau u(x)=\dtau u(x) := \nabla u(x) - \xox \, \dr u ,
\label{Eq:N3}
\eeq
and for a function $n(\omega) \in C^1(S^{d-1})$, we shall consider
$$ 
\nabla_\omega n(\omega) =\domega n (\omega) :=|x|\dtau n(\xox),\qquad \omega=\xox.
$$ 

Let us consider the following assumptions:
\beq
     n =n_1+n_2 \quad\text{with }n_2 \in  L^{\infty},\qquad  n >0 ,
     \label{Eq:A1}
\eeq
\beq
\|(n_1)^{1/2}u\|_2<(1-c_0)\|\nabla u\|_2\quad\text{for all smooth functions $u$ and some $c_0>0$},
     \label{Eq:A1'}
\eeq
\beq
2 \sum_{j\in \ZZ} \; \sup_{C(j)}
\frac{(x\cdot \nabla  n(x))_{-}}{ n(x)}:=\beta < 1.
\label{Eq:A2}
\eeq
Above $(a)_-$ denotes the negative part of $a\in \R.$ 

In \cite{PV} we prove the following result.
%%%%%%%%%%%%%%%%%%%%%%%%%%%%%%%%%%%
\begin{Tm} We assume one of the following two conditions:

\noindent (i) $d\ge 3$, \fer{Eq:A1}--\fer{Eq:A2} and
$R_0=0$;

\noindent (ii) $d= 2$, \fer{Eq:A1}--\fer{Eq:A2}, $n>n_0>0$ and
$R_0= n_0^{-1/2}$;

\noindent Then the solution to the Helmholtz equation \fer{Eq:H}
satisfies,
$$
M^2:= \tn \nabla u \tn ^2_{R_0} + \tn n^{\frac{1}{2}} u
\tn ^2_{R_0} + \int_{|x|\geq R_0}
\f{|\nabla_\tau u|^2}{|x|}
 \, dx\qquad \qquad \qquad
$$
\beq \qquad \qquad \qquad \qquad \qquad \le C (\e +\| n_2\|_{\infty}) \;
N_{R_0}\lp \frac{f}{n^{1/2}}\rp^{2} .
\label{MCEst}
\eeq
\label{Thm1}
\end{Tm}
%%%%%%%%%%%%%%%%%%%%%%%%%%%%%%%%%%%

The homogeneity of the above estimate makes it compatible with the high
frequencies (replace $n$ by $\mu^2 n$).  Our main interest in
\cite{PV} was to obtain estimates with the right scaling.
In particular we were able to recover the well known inequality  for $u$  of Agmon and
H\"ormander in
\cite{AH} in the constant coefficient case. Similar results but not
scaling invariant  were obtained in \cite{JP} and \cite{Zh3}. The scaling
plays a fundamental role in the applications to nonlinear Schr\"odinger
equations (\cite{KPV}) and in the high frequency limit for Helmholtz equations
(\cite{BCKP}, \cite{CPR}).
In this paper we get an explicit Sommerfeld radiation
condition for solutions obtained from the limiting absorption principle. As we have already said
  the estimate of the tangential component of the gradient 
  \beq
  \int_{|x|\geq R_0}
\f{|\nabla_\tau u(x)|^2}{|x|} \, dx<\infty,
\label{tang}
\eeq
given in \fer{MCEst} turns out
to be fundamental. In order to get it we need the smallness assumption given in
\fer{Eq:A2}. We do not know if it is necessary or not. However the condition in \fer{Eq:A2} is necessary and can not  be relaxed to a Coulomb type of decay, even if smallness is added. This is proved in the appendix using as counterexamples a family  of  wave guides for which the estimate of the tangential derivative is false.

In order to prove the new energy estimate, we need to impose some extra assumptions on $n$. 
They are the following ones:
\beq 
\text{there exists}\qquad n_{\infty}(\xox)\in  C^3(S^{d-1}), \quad n_{\infty}(\xox) \geq n_0  >0,
\label{Eq:A3}
\eeq
%%%%%%%%%%%%%%%%%%
and
\beq
|n(x) -n_{\infty}(\xox)| \leq n_\infty(\xox) \,
\f{  \Gamma }{| x|}, \qquad  \Gamma >0\quad n>0 \ ;
\label{Eq:A4}
\eeq
In fact and as we shall prove in section 2 this last assumption \fer{Eq:A4} can 
be largely relaxed if for example $n-n_\infty$ is radial -see \fer {Eq:AA7}, \fer{Eq:A8} below.

We may now state our basic new estimate (already announced in \cite{PVcras}). Its interest relies of course on the 
bounds stated in Theorem \ref{Thm1}.
%%%%%%%%%%%%%%%%%%%%%%%%%%%%%%%%%%%
\begin{Tm}
For dimensions $d\geq 2$, we assume \fer{Eq:A2}, \fer{Eq:A3} and \fer{Eq:A4} and
 use the notation of Theorem
\ref{Thm1}. Then the solution to the Helmholtz equation \fer{Eq:H} satisfies, for $R\geq
R_0$ and $R$ large enough
\beq
\bar M:=\int_{|x|\geq R} |\nabla_\omega n_\infty(\xox) |^2 \; 
\f{|u(x)|^2}{|x|} \, dx \leq
C \ \lc (\e +\| n\|_{\infty})
\; N_{R_0}\lp \frac{f}{n^{1/2}}\rp \, M +M^2 \rc
\label{AngEst}
\eeq
for some constant $C$ independent of $\epsilon$. 
\label{Thm2}
\end{Tm}
%%%%%%%%%%%%%%%%%%%%%%%%%%%%%%%%%%%

We would like to point out the sharpness of this inequality. It says
that the points where $|\nabla_\omega n_\infty(\xox) |$ vanishes on the
sphere are the concentration directions for the energy $|u |^2$. 
   Indeed, we can derive from the
Sommerfeld condition below the following proposition. The proof can
be found at the end
of section 3.

%%%%%%%%%%%%%%%%%%%%%%%%%%%%%%%%%%%
\begin{Prop}
With the assumptions of Theorem \ref{Thm4} below, we have
\beq
{\rm lim_{R\rightarrow \infty}} \f{1}{R} \int_{|x|\leq R} n^{1/2}_\infty(\xox) |u(x)|^2
\, dx =- {\cal I}m \int_{R^d} f(x) \ \bu(x) \ dx.
\label{Eq:Energy}
\eeq
\label{Prop1}
\end{Prop}
%%%%%%%%%%%%%%%%%%%%%%%%%%%%%%%%%%%

Therefore if $u$ vanishes in a neighborhood of the critical points of $n_\infty$ we get the bound
$$
\int_{|x|\geq R} \f{|u(x)|^2}{|x|} \, dx < \infty,
$$
and deduce that ${\cal I}m \int_{\R^d} f(x) \ \bu(x) \ dx =0$. From the Sommerfeld radiation condition and  in the
constant coefficient case this leads (\cite{Ho1}, p. 242) to a
restrictive
condition on the
Fourier transform of $f$, namely
$ \wh f(\xi)=0$ on the sphere $|\xi|= n^{1/2}$. It would be interesting to explore which should be the conclusion in our setting. A natural hypothesis is to substitute  $ \wh f(\xi)$ by the generalized Fourier transform defined by S.
Agmon, J. Cruz-Sampedro and I. Herbst in \cite{ACH} which depends upon
the construction of solutions to the associated eikonal equation.
 The role played by the critical points of $n_\infty$ was
already pointed out by I. Herbst in \cite{He}. There are two related estimates that have been deduced by the same method, the case of stationary kinetic equations (the high frequency limit of Helmholtz equations), see \cite{PVkin}, and the 
case of two layers, see \cite{F}, where the gradient of $n_\infty$ gives a surface Dirac mass. 

Our next purpose is to give an explicit Sommerfeld radiation condition for the solution 
obtained by the limiting absorption principle.  Our result  complements that of Saito in \cite{S}. 
Therefore we shall assume 
\beq
n(x)=\lambda+ p(x),\qquad \lambda>0,
\label{Eq:A7}
\eeq
with $p$  a bounded real function which belongs to $  C^2(\R^n \backslash \{0\})$ and such that
\beq
|\partial^\alpha_x p(x)|\leq c|x|^{-|\alpha|}, \qquad |\alpha|\leq2.
\label{Eq:A6}
\eeq

%%%%%%%%%%%%%%%%%%%%%%%%%%%%%%%%%%%
\begin{Tm}
For dimensions $d\ge 2$, assume  \fer{Eq:A2},  \fer{Eq:A6} and \fer{Eq:A7}. Then
 for $\lambda$ large enough compared to $\|p\|_{{\cal C}^2(\,|x|>1)}$, there
exists a unique solution to the Helmholtz equation with $\e =0$,
$M<\infty,$  and
which satisfies for  any $a>1$,
\beq
\int_{\R^d} \big| \nabla u(x) -i n^{1/2}(x) 
u(x) \ \xox\, \big|^2
\ \f{dx}{1+|x |} \leq C_{a} \int_{\R^d} |f(x)|^2(1+|x|)^a.
\label{1.18}
\eeq
Moreover if there are $n_\infty$, $\Gamma>0$ and $\delta>0$ such that
\beq \label{eq:ndecay}
|n(x) -n_{\infty}(\xox)| \leq n(x) 
\f{\Gamma }{| x|^\delta} \quad \hbox{for $| x|$ large enough},
\eeq
then from  \fer {MCEst} and \fer {1.18} we get
\beq
 \int_{\R^d} \big| \nabla u(x) -i n_\infty^{1/2}(x) 
u(x) \ \xox\, \big|^2
\ \f{dx}{1+|x |} \leq C_{a} \int_{\R^d} |f(x)|^2(1+|x|)^a.
\label{1.188}\eeq
\label{Thm4}
\end{Tm}
%%%%%%%%%%%%%%%%%%%%%%%%%%%%%%%%%%%

Let us compare the above theorem with previously known results. There is a
very extensive literature on the limiting absorption principle, see for
example \cite{E1}, \cite{E2}, \cite{A},  \cite{K}, \cite{BA}, \cite{JP},
\cite{Zh2} and references there in. The situation for the Sommerfeld
radiation condition is different. When
$n=\lambda + V(x)$ and $V$ is a short range potential the question was
settled by Ikebe and Saito in \cite{IS}. Mochizuku and Uchiyama study in
\cite{MU} large range potentials with mild radial oscillations at infinity
like
$V(x)\sim {\rm sin} ({\rm ln} |x|)$. H\"ormander in \cite{Ho2},
chapter XXX, characterizes
the incoming/outgoing solutions  obtained from the limiting absorption
principle by some asymptotic behavior, but in his case
$n_\infty=\lambda=constant.$ More general long range potentials were
considered by Saito in
\cite{S}. Although in this latter work perturbations of first order terms
("magnetic potentials") are also considered let us fix the attention in the
conditions for $V$. Saito writes $V= p+Q$ where $Q$ is a short range
perturbation, while $p$ satisfies \eqref{Eq:A6}. Then he proves a Sommerfeld radiation condition for $\lambda$ large enough
given by $\nabla u\pm i\sqrt\lambda (\nabla \varphi) \ u$, where $\varphi$ is an
appropriate solution for $|x|>R_0$, $R_0$ large enough, of the associated eikonal equation {\footnote
{The existence of this solution was established later on by Barles \cite{Br}.}
\beq|\nabla \varphi|^2=1+\f{p(x)}{\lambda}.
\label{1.20}
\eeq
 Therefore  one cannot expect that in general the vector $\nabla \varphi$ points
at the direction $x/|x|$. An illustrative example is to consider
$$p(x)= - \f{x_1}{  |x|}.$$
In this case and for $\lambda$ large enough, see Remark 1.3 in \cite{S},
$\varphi(x)=a(\lambda)|x|-b(\lambda)x_1$ with $a(\lambda)=1/2[(1+1/\lambda)^{1/2}
+(1-1/\lambda)^{1/2}]$
and
$b(\lambda)=1/2[(1+1/\lambda)^{1/2}
-(1-1/\lambda)^{1/2}]$. This boundary condition differs from ours in all
points except when $\nabla n=0$ (here $n=\lambda+p$). Then the apparent
contradiction is clarified thanks to the estimate \fer{AngEst} which applies for
 this example. 

Notice however that the assumptions for Theorem \ref{Thm2} and for Theorem \ref{Thm4} are different and not comparable.
In the particular case $n=n_\infty$  and regular, no smallness assumption is needed in Theorem \ref{Thm2} because \fer{Eq:A2} is trivially fulfilled, while $\lambda$ has to be large to construct the solution of the eikonal equation, which is a fundamental step in order to prove Theorem \ref{Thm4}.
On the other hand Saito's assumption \eqref{Eq:A6} with $n=\lambda +p$ doesn't imply the existence of the limit $n_\infty$. 

In the proof of Theorem \ref{Thm4} is crucial the estimate of the tangential part of the gradient given in 
\fer{MCEst} to conclude that for $\varphi$ given in \eqref{1.20}
\beq
\int |\nabla_\tau \varphi(x)\, u(x)|^2\,\frac{dx}{1+|x|}<+\infty.
\label{1.21}
\eeq
This is an energy estimate in itself which says that $u$ concentrates along the 
critical points of $\nabla_\tau \varphi$. In section 3 we prove that under some 
conditions these critical points coincide with those of $\nabla_\tau n_\infty$ 
establishing a relation between the energy estimate given in Theorem \ref{Thm2}
 and that in \fer{1.21}. Notice however that  for the proof of Theorem \ref{Thm2} we do not need the existence of a solution to the  eikonal equation \fer{1.20} which well could not exist.

The paper is organized as follows. In section 2 we prove Theorem \ref{Thm2}. Section 3 is devoted to the study of the eikonal equation following \cite{Br}. In particular we give some 
properties of the corresponding solution. The proof of Theorem \ref{Thm4} is
 given in section 4. In Appendix 1 we gather some elementary  identities and in 
 Appendix 2 we give the example of the family of wave guides mentioned above.
\\
\\
\bf Acknowledgments. \rm We want to thank T. Hoffmann-Ostenhoff,  and G.  
Barles for enlightening conversations. We also thank E. Fouassier for pointing out 
the shortcomings of a previous version.

%
%
%%%%%%%%%%%%%%%%%%%%%%%%%%%%%%%%%%%%%%%%%%%%%%%%%%%
\section{Proof of Theorem \ref{Thm2}}
\setcounter{equation}{0}
%%%%%%%%%%%%%%%%%%%%%%%%%%%%%%%%%%%%%%%%%%%%%%%%%%%
%%%%%%%%%%%%%%%%%%%%%%%%%%%%%%%%%%%%%%%%%%%%%%%%%%%
%
%
We are going to prove a  more general version. In fact we will consider two different 
ways of measuring $|n-n_\infty|$. First recall the  assumption \fer{Eq:A4}
\beq
|n(x) -n_{\infty}(\xox)| \leq n(x) \,
\f{  \Gamma }{| x|}, \qquad  \Gamma >0.
\label{Eq:AA4}
\eeq
We can instead assume the alternative conditions
\beq 
|n(x) -n_{\infty}(\xox)| \leq n
\f{\Gamma }{| x|^\delta} \quad \hbox{for $| x|>R_0$}, \qquad \Gamma>0, \quad \delta>0,\quad \hbox{and $R_0>1$};
\label{Eq:AA7}
\eeq

\beq \left\{ \begin{array}{ll}
\hbox{
there exists $\t \beta<1$,  $\delta>0$ and $\bar \Gamma>0$ such that} \\
\lp {|x|\nabla_\tau (n-n_\infty) }\cdot
\f{\partial n_\infty}{\partial \omega}\rp_- \leq
    \ \tilde\beta \, |\f{\partial n_\infty}{\partial\omega}|^2 + n(x)
\f{ \bar \Gamma }{|x|^{\delta}}.
\end{array}
\right.
\label{Eq:A8}
\eeq
In particular if $n-n_\infty$ is radial \fer{Eq:AA7} is sufficient. Also note that from \eqref{Eq:A3} and either \eqref{Eq:AA4}
or \eqref{Eq:AA7} we conclude that 
\beq |n|\leq C\qquad n\geq\frac{n_0}{2}\qquad\hbox{for $|x|$ large enough.}
\label{Eq:AA28}
\eeq

Let us start recalling the fundamental ingredients  of the proof of Theorem 1.1  given in \cite{PV}. In fact in that paper we discarded the estimate for the tangential derivative because it was of not use. However, and as we already said in the introduction, \eqref{tang} turns out to be fundamental in the proof of Theorem 1.2 and Theorem 1.4.

Consider for $R>0$ the special functions
$\Psi,\; \phi$  given
by
\beq \nabla\Psi =  \left\{
\begin{array}{l}
 \f{x}{R}\quad \hbox{for}\quad \vert x \vert \le R,\\
 \; \frac{x}{\vert x \vert} \quad \hbox{for}\quad\vert x \vert \ge R,
\end{array}\right.\label{21}
\eeq
\beq \phi =  \left\{
\begin{array}{l}
 \f{1}{2R}\quad \hbox{for}\quad \vert x \vert \le R,\\
 0 \;\;\; \quad \hbox{for}\quad \vert x \vert \ge R.
\end{array}\right.\label{22}
\eeq
We will also need the following formulas which hold in the distributional sense
\beq D^2_{ij}\Psi=  \left\{
\begin{array}{l}
\f{ \delta_{ij}}{R} \quad \hbox{for}\quad \vert x \vert \le R,\\
 \;\(\delta_{ij}\vert x \vert^{2}-x_{i}x_{j}\)/ \vert x \vert^{3} \quad
 \hbox{for}\quad \vert x \vert \ge R,
\end{array}\right.\label{23}
\eeq
\beq
 \Delta \Psi=  \left\{
\begin{array}{l}
 \f{d}{R} \quad \hbox{for}\quad \vert x \vert \le R,\\
 \;(d-1)/ \vert x \vert \quad
 \hbox{for}\quad \vert x \vert \ge R.
\end{array}\right.\label{24}
\eeq
A simple calculation gives for $d>2$
\beq \frac{1}{4}
\intrd v \; \Delta\(2 \phi-\Delta \Psi\) \ge \frac{d-1}{4 R^2}
\int_{S(R)}v\; d\sigma_{R}\;.\label{25}
\eeq
As it is well known   $\Delta^2\Psi$ is positive for $d=2,$ and has to be treated in a different way. We refer to section 5 in \cite{PV} for the corresponding modifications. 

Adding the identity \eqref{Eq1} to \eqref{Eq3} with the above choices of $\Psi$ and $\phi$ we get 
$$\frac{1}{2R}\intbr \vert \nabla u \vert^2  +  \intbrc\(\vert \nabla
u\vert ^{2} -\vert \frac{x}{\vert x \vert} \cdot \nabla u\vert
^{2}\)/\vert x \vert
$$
$$ + \frac{d-1}{4 R^2} \int_{S(R)}\vert u \vert^2\; d\sigma_{R}+
\frac{1}{2R}\intbr \( n(x) + x\cdot \nabla n(x)\) \vert u \vert^2
$$
$$+ \frac{1}{2}\intbrc \frac{x}{\vert x \vert} \cdot \nabla n(x) \vert
u \vert^2 \le
$$
$$-{\cal R}e\[\f{1}{R}\intbr f(x) \(x\cdot \nabla \bu +\frac{d-1}{2} \bu \)
+  \intbrc f(x)\(\frac{x}{\vert x \vert} \cdot
\nabla \bu + \frac{d-1}{2\vert x \vert} \bu \)\]
$$
\beq
-\e {\cal I}m \[\f{1}{R} \intbr x \cdot \nabla  u \; \bu  +
 \intbrc \frac{x}{\vert x \vert} \cdot \nabla u \; \bu \].
\label{26}
\eeq

We shall estimate separately the various terms in the right-hand side
of the above inequality. We begin by the four terms containing $f$. The pairs containing $f$ and
$\nabla u$ is easy to bound by the right-hand side of \fer{MCEst} because the two norms defined in \eqref{Eq:N1} and \eqref{Eq:N2} 
 are one dual of the other one (see also the next argument). So we concentrate on the other two terms we have for $\delta>0$
 $$\intrd \vert f(x)\vert \;\frac{\vert u\vert}{\vert x \vert}
\le \sum_{j\in \ZZ}\(2^{-j} \int_{C(j)} \frac{ \vert
u\vert^{2}}{\vert x\vert ^{2}}\)^{1/2} \(2^{j} \int_{C(j)} \vert f(x)\vert^{2}
\)^{1/2} \quad \quad
$$
$$\quad \quad \quad \quad \le \( \sup_{R} \frac{1}{R^{2}} \int_{S(R)}\vert
u \vert^2
d\sigma_{R} \)^{1/2} \sum_{j\in \ZZ} \(2^j \int_{C(j)} \vert f(x)\vert^{2}
\)^{1/2}$$
\beq
 \quad \quad \le \delta \sup_{R} \frac{1}{R^{2}} \int_{S(R)}\vert u \vert^2
d\sigma_{R} + C_{\delta} N(f)^{2}.
\label{266}
\eeq

We consider now the terms on
$n$.  We have,
$$ \frac{1}{R} \intrd (\nabla \Psi \cdot \nabla  n)_{-}\vert u \vert^2
\le \sum_{j\in \ZZ} \int_{C(j)} n \vert u \vert ^{2} \;
\frac{(x \cdot \nabla  n)_{-}}{\vert x \vert n}
$$
\beq
\le \( \sup_R\frac{1}{R} \int_{B(R)} n \vert u \vert^2 \)
\sum_{j \in \ZZ} \sup_{C(j)}2^{j+1}
\frac{(x \cdot \nabla  n)_{-}}{\vert x \vert n}
\le \beta_{1} \sup_{R} \frac{1}{R} \intbr n \vert u \vert^2.
\label {27}
\eeq

Finally we have to bound the terms involving $\e$.
 From the Helmholtz
equation one deduces
$$ \e \intrd  \vert u \vert^2 \le \intrd \vert f\; \bu \vert,$$
$$ \intrd \vert \nabla u\vert^{2} \le
\intrd n \vert u \vert^2 + \intrd \vert f\; \bu \vert.$$
With the assumptions  \eqref{Eq:A1} and  \eqref{Eq:A1'} we obtain
$$\intrd \vert \nabla u\vert^{2} \le C \( \intrd n_{2} \vert u \vert^2 +
\intrd
\vert f\; \bu \vert \).$$
Using Cauchy-Schwarz inequality and the above
inequalities we get
$$  \e \intrd  \vert \nabla u\vert \; \vert u \vert \le C \e^{1/2}
\( \intrd n_{2}\vert u \vert^2+\intrd  \vert f\; \bu \vert\)^{1/2}
\(\intrd  \vert f\; \bu \vert \)^{1/2}
$$
$$\le C \e^{1/2} \intrd \vert f\; \bu \vert + C \(\e \sup \vert
n_{2}\vert \; \intrd \vert f\; \bu \vert  \; \intrd  \vert u \vert^2 \)^{1/2}
$$
$$\le C (\e + \sup \vert n_{2}\vert)^{1/2} \intrd \vert f\; \bu \vert
$$
$$\le C (\e + \sup \vert n_{2}\vert)^{1/2} \tn n^{1/2} u\tn
N(\frac{f}{n^{1/2}})
$$
\beq
\le \delta \tn n^{1/2} u\tn ^{2} +C_{\delta}\; (\e + \sup \vert
n_{2}\vert)\;
(N(\frac{f}{n^{1/2}}))^2.
\label{28}
\eeq
Then plug \eqref{266}-\eqref{28} into \eqref{25} and take the supremum in $R$. Then we get
\eqref{MCEst} and in particular \eqref{tang}.
\\

 Let us prove now Theorem 1.2. The proof consists in using the basic equality \fer{Eq3}  with a test function that depends on the
behavior of $n(x)$ at infinity. We choose for $R\geq R_0$ such that \eqref{Eq:AA28} holds and define
$$\Psi_q(x)= q(\f{|x|}{R}) \ n_\infty(\xox)
$$
for some non-decreasing smooth  function $q(r)=0$ for $r\leq
1$ and $q(r)=r$ for $r\geq 2$.

With this choice, we will show that the only new information (compared to Theorem
\ref{Thm1}) in \fer{Eq3}, is given by the term
\beq
\int_{R^d} \nabla n(x)\cdot \nabla \Psi_q (x) \ | u(x)|^{2} \, dx.
\label{Eq:AA29}
\eeq
Then  we will take care of the other terms using \fer{MCEst}.

As a {\it first step} we consider \eqref{Eq:AA29}. We simplify the notation using
$q=q(\f{|x|}{R})$. Then we get
$$
\int_{R^d} \nabla n(x)\cdot \nabla \Psi_q (x) \ | u(x)|^{2} \, dx
=\int_{R^d} q \ |\domega n_\infty(\omega)|^2 \ \f{| u(x)|^{2}}{|x|^2} \, dx
$$
\beq
+ \int_{R^d} \dr n(x) \ \f{q'}{R} \ n_\infty(\omega) \ | u(x)|^{2} \, dx
+ \int_{R^d} q |x|\nabla_\tau\lp n(x)-n_\infty(\xox) \rp \  \domega n_\infty(\omega)
\ \f{| u(x)|^{2}}{|x|^2} \, dx.
\label{2.1}
\eeq
The first term on the right-hand side gives the control we look for.
As for the second term we just have to consider the negative part of $ \dr n(x) $. Then we proceed as
in \eqref{27} to  get the lower bound
$$
- \f{C}{R} \ \beta \ \| n_\infty \|_{L^\infty} \ \tn n^{1/2} u \tn_{R_0}^2.
$$
Let us consider first condition \fer{Eq:AA4}. Then the last term in \fer{2.1}
is bounded below as follows. After integration by parts, it is also
given by
$$
- {\cal R}e \int_{R^d}  \f{ q}{|x|^2} \ \lp n(x)-n_\infty(\xox) \rp
\lp  D^2_\omega  n_\infty \ | u(x)|^2 +
2 \ \domega n_\infty (\omega)\ |x|\nabla_\tau u(x) \ \bu \rp \, dx
$$
$$\geq -  \Gamma \lc \| n_\infty \|_{C^2} \int_{\R^d} q \ n(x) \ \f{|
u(x)|^2}{|x|^3} \,
dx + 2\| n\|_{L^\infty} \int_{\R^d} \f{ q }{|x|^2} \ |\domega n_\infty (\omega)|
\ |\nabla_\tau
u(x)| \ |u | \, dx \rc
$$
$$ \geq - \f{ C}{R^2} \ \| n_\infty \|_{C^2} \tn n^{1/2} u\tn_R^2
- \f12 \int_{\R^d}q \ |\domega n_\infty(\omega)|^2 \ \f{|
u(x)|^{2}}{|x|^2} \, dx
$$
$$ - \ \f{ C}{R} \ \| n \|^2_{L^\infty}
\int_{|x|\geq R}
\f{|\nabla_\tau  u(x)|^2}{|x|} \, dx.
$$
%%%%%%%%
Let us assume now \fer{Eq:AA7}, then the last term of \fer{2.1} is
bounded below by
$$- \t \beta \int_{\R^d} q \ |\domega n_\infty(\omega)|^2 \ \f{|
u(x)|^{2}}{|x|^2} \, dx
-  \f{C}{R^\delta} \tn n^{1/2}u\tn_{R_0}^2.
$$

As a conclusion of this first step we have obtained
$$
\int_{\R^d} q(\f{|x|}{R}) \ |\domega n_\infty(\omega)|^2 \ \f{| u(x)|^2}{|x|^2} \,
 dx
$$
\beq
\leq C_1 \int_{\R^d} \nabla n(x)\cdot \nabla \Psi_q (x) \ | u(x)|^{2} \, dx
+C_2  \ M^2.
\eeq

\medskip

The {\it second step} is to provide a control on all the remaining terms thanks to the basic equality \fer{Eq3}.

We have
$$
\nabla \bu(x) \cdot D^{2}\Psi_q (x) \cdot \nabla u(x)=
\f{q'}{R|x|} \  n_\infty(\omega) \  |\nabla_\tau u|^2 +
\f{q''}{R^2} \  n_\infty (\omega) \  |\dr u|^2
$$
\beq
+2 \ (\f{q'}{R|x|}-\f{q}{|x|^2}) \ {\cal R}e [ \dr \bu \
\domega n_\infty(\omega)  \cdot \nabla_\tau u]
+ \f{q}{|x|^2} \ \nabla_\tau \bu \cdot D^2_\omega n_\infty (\omega) \cdot \nabla_\tau u.
\label{2.3}
\eeq
And because the terms $(\f{q'}{|x|}-\f{q}{|x|^2})$ and $q''$ are
supported in the ball $\{|x| \leq R\}$, we see that from \eqref{MCEst}  all the terms in the
corresponding integral are
bounded  by $C \|n_\infty\|_{C^2} M/R$.

Next, we consider the term
$$\int_{\R^d} \Delta^{2}\Psi_q(x) | u(x)|^{2}=-2{\cal R}e \int_{\R^d} \nabla
\Delta \Psi_q \cdot \nabla u \, \bu
$$
\beq
\leq 2 \int_{\R^d} \f{|\nabla \Delta \Psi_q|}{n^{1/2}}\, |\nabla u| \, |n^{1/2}u|.
\label{Eq:230}
\eeq

But we have
$$ |\nabla  \Delta \Psi_q(x) | \leq C \ \| n_\infty \|_{C^3} [\f{q}{|x|^3}
+ \f{q'}{R |x|^2}]\leq c\frac{\|n\|_{\infty}}{|x|^2}.
$$
Therefore from \eqref{MCEst} we get that \eqref{Eq:230} is bounded. by
$$C \|\f{n_\infty}{n_0^{1/2}} \|_{C^3}  M/R.$$

The right-hand side terms containing $f$  can be treated as in \eqref{266} and are respectively
upper bounded  by 
$$C \| n_\infty \|_{C^1} N_{R_0}(f) \tn \nabla u \tn_{R_0},$$
and by
$$C \|\f{n_\infty}{n_0^{1/2}} \|_{C^2} N_{R_0}(f) \tn
n^{1/2}u \tn_{R_0}.$$

The last term to be bounded is 
\beq
\e {\cal I}m\int_{\R^d}  \nabla \Psi_q (x) \cdot\nabla \bu (x) u(x).
\label{JFA1}
\eeq
But this is done as in \eqref{28}.
The proof of Theorem \ref{Thm2} is over.

%%%%%%%%%%%%%%%%%%%%%%%%%%%%%%%%%%%%%%%%%%%%
\section{The eikonal equation}
\setcounter{equation}{0}
\label{sec:ccc}
%%%%%%%%%%%%%%%%%%%%%%%%%%%%%%%%%%%%%%%%%%%%

In order to determine the phase arising in the Sommerfeld radiation
 condition, we need to solve the eikonal equation
\beq
\left\{\begin{array}{l}
|\nabla \vp |^2 =\f{1}{\lambda} n(x)= 1+ \f{1}{\lambda} p(x), \qquad x \in \R{^d},
\quad \lambda>0,
\\
\\
\vp(0)=0, \qquad \nabla \vp(x)/|x| \;  \text{ bounded},
\end{array} \right.
\label{eq:eik} 
\eeq
where we assume that $\lambda>0$ is large enough and that $p \in C^2(\R{^d})$ and satisfies 
\beq
|D^\al p(x)| \leq C (1+ |x|)^{-\al}, \qquad \text{for } |\al| \leq 2.
\label{as:hj1}
\eeq

In order to take into account the linear growth at infinity (which is essential for uniqueness), the  unique viscosity solution $\vp$ to \fer{eq:eik} is better described in terms
 of the  bounded function  $g(x)$ defined as
$$
\vp(x)= |x| g(x) .
$$ 
This change of unknown yields the following Hamilton-Jacobi equation for $g$:
\beq
\left\{\begin{array}{l}
|g|^2 + 2r g \p_r g + |x|^2 |\nabla g|^2= 1+ \f{1}{\lambda}p(x), \qquad x \in \R
{^d}\backslash \{0\},
\\
g(0)=\sqrt{\f{n(0)}{\lambda}}.
\end{array} \right.
\label{eq:hj}
\eeq
From this H.-J. equation, one can derive properties of the phase $\vp$, which we state now.

%--------------------------------------------------
\subsection{Statements of the results}
%-------------------------------------------------

Let us recall the results of  Barles given in \cite{Br}.

\begin{theorem} [Barles, \cite{Br}]  Under assumption \fer{as:hj1}, there exists a unique  solution to 
\fer{eq:hj}. Moreover it satisfies the following estimates, for $\lambda$ large enough and $x
 \neq 0$,
$$
\sqrt{ 1-  \min p/\lambda} \leq g(x) \leq  \sqrt{ 1+  \max p/\lambda},
$$
\beq
|x|  |D  g(x)|  \leq C/ \lambda, \quad |x|^2  |D^2 g(x)|  \leq C/\lambda, \quad 
|x|^3  |D^3 g(x)|  \leq C/\lambda.
\label {3.11}
\eeq
\end{theorem}

\medskip

In the proof of Theorem \ref{Thm4} we shall derive  the estimate
\beq
\int |\nabla_\tau\varphi\, u|^2\frac{dx}{1+|x|}<+\infty,
\label{311}
\eeq
which is similar to our estimate \fer{AngEst}  in Theorem \ref{Thm2}. Therefore it is a natural question to know which is the relation between the two inequalities \eqref{AngEst} and
\eqref{311}. In fact we shall see that under some extra assumptions on the index $n$, both inequalities are equivalent.
Indeed, assume
\beq
| \p_r p(x)| \leq C r^{-1-\delta}, \qquad \text{for some } \delta \in ]0,1[ ,
\label{as:hj2}
\eeq
with $r=|x|$ and $\p_r= \f{x}{|x|} \cdot \nabla$ the radial derivative.
Firstly, as a consequence of this assumption, $n(x)$ admits a radial limit
\beq
n(r \f{x}{|x|}) \to n_\infty (\f{x}{|x|})\quad  \text{as }\;  
r \to \infty,
\label{eq:radlimt}
\eeq
and in fact the decay rate in our previous assumption \fer{eq:ndecay} also follows from \fer{as:hj2}. Then, we have
%---------------------------
\begin{theorem} Under assumptions \fer{as:hj1} and \fer{as:hj2}, the solution to
\fer{eq:hj} satisfies for $\lambda$ large enough and $x \neq 0$ the  estimate, 
\beq
r^{1+\delta} |\p_r g(x)| \leq C/ \lambda,  
\label{eq:hjest}
\eeq 
and, as in \fer{eq:radlimt}, $g(r \f{x}{|x|}) \to g_\infty (\f{x}{|x|})$ as $r \to \infty$, a smooth solution to the equation 
\beq
g_\infty(\omega)^2 + |\nabla_\omega g_\infty(\omega)|^2 =n_\infty(\omega), \qquad \omega \in \Sph^{d-1}.
\label{eq:hjlim}
\eeq
Moreover
$$
|\nabla_\tau\varphi|=|\nabla_\omega g_\infty(\omega)|+ O( r^{-\delta}).
$$
\label{th:radlim}
\end{theorem}

\noindent {\bf Remark.} The solution to the limiting equation \fer{eq:hjlim} enjoys a specific property.  
 From \fer{3.11} we get that $(g_\infty I+D^2_\omega g_\infty)$ is invertible and  
 after differentiation in \fer{eq:hjlim} we obtain
 \beq 
0< c_1|\nabla_\omega g_\infty| \leq |\nabla_\omega n_\infty|  \leq c_2|\nabla_\omega g_\infty|.
 \label{3.13}
\eeq
Therefore the critical points (and thus the extrema) of $n_\infty$ and $g_\infty$ coincide.
\\

The derivation of Theorem \ref{th:radlim} can be seen either from the  representation 
through the method of characteristics or from the more modern PDE point of view. We give both proofs  for the sake of completeness.
%---------------------------------------
\subsection{Proof of Theorem \ref{th:radlim} based on bicharacteristics}
%---------------------------------------

We decompose the proof in two steps. We first recall the definition of the characteristics and how they give a representation formula for the solution to   the eikonal equation \fer{eq:eik}, then we come to the proof of the estimates in  Theorem  \ref{th:radlim}.
\\

\noindent {\em First step. Characteristics}. The bicharacteristics are classically defined as the solutions to the differential system, parameterized by some $q\in \Sph^d$, 
\beq \label{sys:bichar}
\left\{ \begin{array}{ll}
\dot X(t;q)= 2 P(t;q), & X(t=0;q)=0,
\\ \\
\dot P(t;q)= \f{1}{\lambda}\nabla n\big( X(t;q) \big), \quad & P(t=0;q)= \f{1}{\lambda^{1/2}}q \sqrt{n(0)}, \; |q|^2=1,  
\\ \\
\dot \Phi(t;q)=  \f{2}{\lambda} n\big(X(t;q)\big), &  \Phi(t=0;q)= 0. 
\end{array} \right.
\eeq
We note for later purpose that 
$$
\f{d}{dt} [|P|^2- \f{1}{\lambda}n(X)]= 0, \; \text{ therefore } \;  |P(t;q)|^2 = \f{1}{\lambda}n\big(X(t;q)\big) \; \forall t \in \R.
$$
Assume first that we can define a unique diffeomorphism (this involves smallness assumptions)
\beq \label{eq:diff}
x\in \R^d\backslash \{0\} \to (t,q) \in ]0, \infty[\times \Sph^{d-1}, \qquad x=X(t;q) .
\eeq
Then, we recall the standard representation formula for the solution $\vp(x)$ to \fer{eq:eik} (see \cite{lions} for details and complements).
%--------------------------------------------
\begin{lemma} For all $t >0$ and $q\in  \Sph^{d-1}$ such that condition \fer{eq:diff} holds, we have 
$$
P(t;q)= \nabla \vp \big(X(t;q)\big), \qquad \vp\big(X(t;q)\big)=  \Phi (t;q),
$$
for some smooth solution $\vp(x)$ to \fer{eq:eik}.
\label{lm:bichar}
\end{lemma}
%---------------------------------------------
\proof  Because both the H.-J. solution and the differential system \fer{sys:bichar} are stable by smooth perturbations of $n$, we can assume that $n$ is constant in a neighborhood of the origin. Then, we divide the proof  in three steps.
\\
(i) Then, close to $x=0$ the solution to the H.-J. equation is given  $\vp(x)=  \f{1}{\lambda^{1/2}}\sqrt{n(0)} \; |x|$ and the bicharacteristics are 
$$
X(t;q)= 2 t\; q \sqrt{n(0)}, \ P= q \sqrt{n(0)}=\nabla \vp(x), \ \Phi= 2 n(0) t= \sqrt{n(0)}|X|=\vp(x).
$$
Hence the result holds true. 
\\
\\
(ii) We prove that $P(t;q)$ is a gradient, i.e., $\Big( \f{\p P}{\p x}\Big)$ is a symmetric matrix. 
Indeed,  from (i), it is true for $t \approx 0$ and it remains to see that for all times 
$\f{d}{dt} \Big( \f{\p P}{\p x}\Big)$ is a symmetric matrix. But we can compute 
$$
 \begin{array}{rl}
\f{d}{dt} \Big( \f{\p P}{\p x}\Big)&= \f{d}{dt}\Big[  \Big( \f{\p P}{\p (t,q)}\Big)\Big( \f{\p X}{\p (t,q)}\Big)^{-1} \Big]
\\
\\
&= \Big( \f{\p \dot P}{\p (t,q)}\Big)\Big( \f{\p X}{\p (t,q)}\Big)^{-1}+  \Big( 
\f{\p P}{\p (t,q)}\Big)  \f{d}{dt} \Big( \f{\p X}{\p (t,q)}\Big)^{-1}
\\
&= \f{1}{\lambda}D^2 n\big(X(t;q)\big) \Big( \f{\p X}{\p (t,q)}\Big) \Big( \f{\p X}{\p (t,q)}\Big)^{-1}
- \Big( \f{\p P}{\p (t,q)}\Big) \Big( \f{\p X}{\p (t,q)}\Big)^{-1} \Big( \f{\p \dot X}{\p (t,q)}\Big) \Big( \f{\p X}{\p (t,q)}\Big)^{-1}
\\
&= \f{1}{\lambda}D^2 n\big(X(t;q)\big) -2 \Big( \f{\p P}{\p (t,q)}\Big) 
\Big( \f{\p X}{\p (t,q)}\Big)^{-1} \Big( \f{\p P}{\p (t,q)}\Big) \Big( \f{\p X}{\p (t,q)}\Big)^{-1} 
\\
&= \f{1}{\lambda}D^2 n\big(X(t;q)\big) - 2\Big( \f{\p P}{\p x}\Big) \Big( \f{\p P}{\p x}\Big)
\\
&=   2P\cdot \Big( \f{\p^2 P}{\p x^2}\Big),
\end{array} 
$$
because $|P|^2=\f{1}{\lambda}n$ implies 
$\Big( \f{\p P}{\p x}\Big) \Big( \f{\p P}{\p x}\Big)+ P\cdot \Big( \f{\p^2 P}{\p x^2}\Big) 
= \f{1}{2\lambda} \Big( \f{\p^2 n}{\p x^2}\Big)$. 
We have obtained indeed that $\f{d}{dt} \Big( \f{\p P}{\p x}\Big)$ 
is a symmetric matrix and our claim is proved. 
\\
(iii) From step (ii), we can write  $P(t;q)= \nabla \vp \big(X(t;q)\big)$for some function $ \vp$. And since 
$|P|^2=\f{1}{\lambda}n$, we obtain that $\vp$ solves the H.-J. equation. It remains to identify $\Phi$ and $\vp$. To do so, we write 
$$
\f{d}{dt} \vp\big(X(t;q)\big) =\nabla \vp\big(X(t;q)\big) \; \f{d}{dt} X(t;q)= 2
 |P(t;q)|^2 = \f{2}{\lambda}n\big(X(t;q)\big).
$$
This means that $\vp\big(X(t;q)\big) = \Phi(t;q)$ since the identity holds at $t=0$ and the two functions have the same derivatives. This completes the proof of Lemma \ref{lm:bichar}.

\bigskip
The reciprocal to Lemma \ref{lm:bichar} is also true. If the H.-J. solution is smooth, then we can derive that  the bicharacteristic system satisfies the invertibility condition. As before, for simplicity we assume that $n$ is constant in a neighborhood of the origin.
%--------------------------------------------
\begin{lemma} Assume that the solution $\vp$ to H.-J. equation is smooth and satisfies $x\cdot \nabla \vp(x) \geq c |x|$ and consider the autonomous differential equation 
\beq
\f{d}{dt} X(t;q)= \nabla \vp (X(t;q)), \qquad X(0;q)=0, \quad X(t;q)\approx n^{1/2}(0) qt \; (t \; small).
\label{eq:clphi}
\eeq
This gives the solution to \fer{sys:bichar} with the formulas of Lemma \ref{lm:bichar}.
 
\end{lemma}
%---------------------------------------------
\proof  
We prove successively that: 
\\
(i) it is the solution to the bicharacteristic system; 
\\
(ii) $(t>0,q \in \Sph^{d-1})\mapsto X(t;q)\in \R^d\backslash \{0\}$ is one-to-one;
\\
(iii) it is a diffeomorphism.
\\

The point (i) follows from the computation in the Lemma \ref{lm:bichar} and we do not reproduce them again.
\\

The  point (ii) can be proved as follows.
\\
(a) into : if $X(t_0;q)=X(t_0+\tau;q')$, by uniqueness of the system for backward time $X(t;q)=X(t+\tau;q')$ for all $t$ and thus (considering what happens close to the origin) $\tau=0$ and $q=q'$.
\\
(b) onto: take a point $y \neq 0$, and solve the equation \fer{eq:clphi} backward from say $t=0$ and call $Y(t,y)$ the solution. Because we know that $x\cdot \nabla \vp(x) \geq c |x|$, we have
$$
\f{d}{dt} \f{|Y(t)|^2}{2}\geq c |Y(t)|.
$$
Therefore it reaches $Y(\tau;y)=0$ in finite time ($\tau <0$). Now, the system \fer{eq:clphi} is autonomous. Therefore we can solve it from $t=0$, with $q_0=n(0)^{-1/2} \lim_{s\to 0}  Y(\tau+s)/s$ and one has $y=X(-\tau; q_0)$. 
\\

For point (iii), we just notice that the above construction, based on the Cauchy-Lipschitz theory, also provides $C^1$ regularity and $C^1$ regularity of the inverse.  
\\

These two lemmas show that, in the regime of assumption \fer{as:hj1}, with $\lambda$ large enough, the solution to the H.-J. equation is given by the bicharacteristics. We can use this fact to conclude the proof of Theorem \ref{th:radlim}.
\\
\\
  
\noindent {\em Second step. Estimates in Theorem \ref{th:radlim}}. We divide the proof of these estimates in two steps. 
First an estimate on the rays is given, then we prove \fer{eq:hjest}.  Then, the 
second statement, \fer{eq:hjlim}, is an easy consequence obtained passing to the limit  as $r\to \infty$.
\\
(i) We have: $|X(t;q)|/t =O(1)$. Indeed, we can compute
$$
\f{d}{dt} \f{|X|^2}{2} = X \cdot \dot X =2  X \cdot P,
$$
and
$$
\f{d}{dt} X \cdot P = 2 |P|^2+ \f{1}{\lambda}X\cdot \nabla n= \f{2}{\lambda} n+ \f{1}{\lambda}X\cdot \nabla n.
$$
As a consequence of assumption \fer{as:hj2}, we have
$$
\f{d^2}{dt^2} \f{|X|^2}{4} = \f{1}{\lambda}n +o(1) \quad \text{for $X$ large},
$$
and thus the result (i) is proved.
\\
\\
(ii) We have:  $|x|^{1+\delta} \p_r g$ is bounded. Using the calculation of step (i), we have 
$$
\begin{array}{rl}
\Phi(t;q)&= \f{2}{\lambda} \int_0^t n\big(X(s;q)\big)ds= \int_0^t [\f{d}{dt} X \cdot P-\f{1}{\lambda}X\cdot \nabla n]ds
\\ \\
&= X \cdot P(t;q)-F(t;q)
\end{array} 
$$
with 
$$
\qquad F(t;q)= \f{1}{\lambda}\int_0^t X(s;q) \cdot \nabla n\big(X(s;q)\big) \;ds, \qquad |F(x)| \leq 
C[1+ |x|]^{1-\delta}. 
$$
In other words, for $x$ large we have
$$
\p_r g = \f{1}{|x|} \p_r \vp- \f{\vp}{|x|^2}=  \f{F(x)}{|x|^2}=O\big( \f{1}{|x|^{1+\delta}}\big),
$$
and the claim (ii) is proved. This concludes the proof of \fer{eq:hjest} and thus the proof of Theorem \ref{th:radlim} is complete. 

%%%%%%%%%%%%%%%%%%%%%%%%%%%%%%%%%%%%%%%%%%%%%%%%%%
%---------------------------------------
\subsection{Proof of Theorem \ref{th:radlim} based on H.-J. equations}
%---------------------------------------

The proof is based on the equation for $\p_r g$,
$$
2 g \p_r g+  r (\p_r g)^2 +  r g \p^2_{rr} g +  r |\nabla g|^2  +   r^2 \nabla g
 \cdot \nabla \p_r g =\f{1}{2\lambda} \p_r n.
$$

\noindent {\it Upper bound}. We consider the maximum point (if it is not attained, 
then perturbation methods, see \cite{lions}, apply)
$$
\max [r^{1+\delta} \p_r g(x) ]= r_0^{1+\delta} \p_r g(x_0), 
$$
and at the point $x_0$, we have
$$
(1+\delta)  \f{x_0}{|x_0|} \p_r g(x_0)+ r \nabla \p_r  g(x_0) =0.
$$ 
Inserting this in the equation on $\p_r g$, we find
$$
2 g \p_r g + r (\p_r g)^2  - (1+\delta) g  \p_r g+ r |\nabla g|^2  -  (1+\delta)
  r \nabla g \cdot  \f{x}{|x|} \p_r g(x_0)  =\f{1}{2\lambda} \p_r n.
$$
Therefore we obtain 
$$
(1- \delta) [ g \p_r g + r (\p_r g)^2] \leq \f{1}{2\lambda} \p_r n,
$$
and thus we arrive at the upper bound 
$$
\max [ r^{1+\delta} \p_r g(x)] \leq  \f{1}{2(1- \delta)\lambda}  \max [ r^{1+\delta} \p_r n(x) ]/\min g.
$$
\\
\noindent {\it Lower bound}. With the same calculation as above, at the point $x_1$ where 
the minimum is attained
$$
\min [r^{1+\delta} \p_r g(x) ]= r_1^{1+\delta} \p_r g(x_1), 
$$ 
we have successively
$$
(1- \delta) g \p_r g - \delta r (\p_r g)^2 +  r |\nabla g|^2  =  \f{1}{2 \lambda} \p_r n,
$$
$$
(1- \delta) g \p_r g - \delta r^{-(1+2\delta)} \big(r^{1+\delta} \p_r g\big)^2 +\f 1 r  [\frac{n}{\lambda}
-g^2-2rg \p_r g] = \f{1}{2 \lambda} \p_r n ,
$$
$$
(1+ \delta) r^{1+\delta} g \p_r g + \delta r^{-\delta} \big(r^{1+\delta} \p_r g\big)^2 
=r^\delta  [\frac{n}{\lambda} -g^2] -  \f{r^{1+\delta}}{2\lambda} \p_r n,
$$
which gives the result.

%
%%%%%%%%%%%%%%%%%%%%%%%%%%%%%%%%%%%%%%%%%%%%%%%%%%%
\section{Proof of Theorem \ref{Thm4} and of
Proposition \ref{Prop1}}
\label{sec:proof}
\setcounter{equation}{0}
%%%%%%%%%%%%%%%%%%%%%%%%%%%%%%%%%%%%%%%%%%%%%%%%%%%
%%%%%%%%%%%%%%%%%%%%%%%%%%%%%%%%%%%%%%%%%%%%%%%%%%%

We will need the following uniqueness theorem.
\begin{Tm}
For dimensions $d\ge 2$, assume \fer{Eq:A1}-\fer{Eq:A2}. If $u$ is a solution of
 $\Delta u+nu=0$ with $u$ and
$\nabla u$ locally in $L^2$, and such that
\beq
{\rm liminf} \int_{|x|=r} \lp | \nabla u|^2+| u |^2\rp \  d\sigma(x) = 0, \qquad {\rm as} \quad r \to \infty,
\label{Eq:0}
\eeq
then $u=0.$
\label{Thm0}
\end{Tm}
\noindent {\bf Proof.} 
Assume $R_j$ is a sequence going to infinity such that
\beq
{\rm lim_j} \int_{|x|=R_j} \lp | \nabla u|^2+| u |^2\rp \  d\sigma(x) = 0, \qquad {\rm as} \quad j \to \infty.
\label{Eq:00}
\eeq
Consider $\Psi$ and $\phi$ as in \eqref{21} and \eqref{22}. Then fix $j$ and use the  multiplier $\nabla\Psi \nabla\bar u+(1/2)\Delta\Psi \bar u+\phi\bar u$  in Helmholtz equation
 $$\Delta u+nu=0,$$
 in the region $|x|\leq R_j$, and for $R\leq R_j$.
 Then we repeat the procedure given at the beginning of section 2 to prove Theorem 1.1. 
 Note that just the multiplier which involves $ \Psi$ will create boundary terms after the integration  by parts. Therefore the basic identity \eqref{Eq3} has to be modified. In this case we shall obtain
 \begin{equation}\label{res0}
\begin{split}
 0=&\int_{B(0,{R_j})} {\nabla }\bar{u}(x)D^2 \Psi(x)
\nabla u(x)dx \\
&-\frac{1}{4}\int_{B(0,{R_j})} |u(x)|^2 {\Delta}^2 \Psi (x)dx +
 \frac{1}{2} \int_{B(0,{R_j})}|u(x)|^2(\nabla n \cdot
\nabla \Psi ) (x)  dx \\
&+ \frac{1 }{2} \int_{S_{R_j}}
\partial_r \Psi (x) |\nabla u(x)|^2 d \sigma (x)   -
   \Re  \int_{S_{R_j}}
\partial_r u(x) \nabla \Psi (x) \cdot \nabla \bar u(x) d \sigma (x)
\\
&+\frac{1 }{4}  \int_{S_{R_j}}
\partial_r \Delta\Psi (x)|u(x)|^2 d \sigma (x)   -
  \frac{1 }{4} \int_{S_{R_j}}
\partial_r (|u(x)|^2) \Delta \Psi (x) d \sigma (x)   \\
&+\frac{1}{2}  \int_{S_{R_j}}
\partial_r \Psi (x)|u(x)|^2 n(x)d \sigma (x).
 \end{split}
\end{equation} 
The inner terms are the same as  in the proof of Theorem 1.1.  Therefore adding those given by $\phi$
we shall obtain following the same argument given in section 2
\begin{equation}
\begin{split}
\sup_{R\leq R_j} \frac{1}{R}\int_{B(0,R)} &({|\nabla }{u}(x)|^2dx +n|u(x)|^2 )dx \\
 &\leq C \int_{S_{R_j}}
|\partial_r \Psi (x) |\,|\nabla u(x)|^2 d \sigma (x)   
+|   \Re  \int_{S_{R_j}}
\partial_r u(x) \nabla \Psi (x) \cdot \nabla \bar u(x) d \sigma (x)|\\
&+C  \int_{S_{R_j}}
|\partial_r \Delta\Psi (x)|\,|u(x)|^2 d \sigma (x)   +
 C \int_{S_{R_j}}
|\partial_r (|u(x)|^2)|\,| \Delta \Psi (x)| d \sigma (x)   \\
&+C  \int_{S_{R_j}}
|\partial_r \Psi (x)|\,|u(x)|^2 n(x)d \sigma (x).
 \end{split}
\end{equation} 
Then taking the limit in $j$ the theorem follows from \eqref{Eq:00}.
\newline
 
%%%%%%%%%%%%%%%%%%%%%%%
%%%%%%%%%%%%%%%%%%%%%%%%%

\noindent {\bf Proof of Theorem \ref{Thm4}.}
Let us recall the estimate  obtained by Saito in Theorem 1.6 of \cite{S}. There is $a$, with 
$1<a\leq2$ such that
\beq
\int\,|\nabla u-i\lambda^{1/2}\nabla\varphi\,u|^2\,\f{dx}{(1+|x|)^{2-a}}\,\leq c\int|f|^2(1+|x|)^a\,dx.
\label{4.10}
\eeq
In order to prove this inequality Saito needed the existence of  $\varphi$, a solution of  the eikonal equation for $|x|>R_0$ with $R_0$ large enough, which was established later on by Barles in \cite{Br}, and that we gave in Theorem 3.1 of section 3. 

From \eqref{4.10} we get
\beq \int |\nabla u -i\lambda^{1/2}\nabla\varphi u|^2 \f{dx}{1+|x|}\leq c\int|f|^2(1+|x|)^a\,dx.
 \label{4.16}
\eeq
Hence looking at just the tangential parts of the above inequality
and from \eqref{tang} we get
$$ \int |\lambda^{1/2}\nabla_\tau\varphi u|^2\,\frac{dx}{1+|x|}\leq c\int|f|^2(1+|x|)^a\,dx.
$$
From the eikonal equation we have that
$$
n-\lambda|\partial_r\varphi |^2=|\lambda^{1/2}\nabla_\tau\varphi |^2.
$$
Recall that from the properties of $\nabla\varphi$ \fer{3.11} we know that 
$\partial_r\varphi=g(x)+O(1/\lambda)>0$.
Then
$$
|\lambda^{1/2}\partial_r\varphi-n^{1/2}|=\frac{|\lambda^{1/2}\nabla_\tau\varphi |^2}{|\lambda^{1/2}\partial_r\varphi+n^{1/2}|}
\leq c |\lambda^{1/2}\nabla_\tau\varphi |^2.
$$
Also looking at the radial part in (\ref{4.16})  we get 
$$\int| \partial_r  u-i\lambda^{1/2}\partial_r\varphi u|^2\,\frac{dx}{1+|x|}\leq c\int|f|^2(1+|x|)^a\,dx.
$$
Finally from the above estimates
we get
$$
\int |\nabla u-in^{1/2}\frac{x}{|x|}u|^2\,\frac{dx}{1+|x|}=
\int(| \partial_r  u-in^{1/2}u|^2+|\nabla_\tau u|^2)\,\frac{dx}{1+|x|}
$$
$$
\leq c\int(| \partial_r  u-i\lambda^{1/2}\partial_r\varphi u|^2+
|(\lambda^{1/2}\partial_r\varphi-n^{1/2})u|^2+ |\nabla_\tau u|^2)\,\frac{dx}{1+|x|}
$$
$$\leq c\int|f|^2(1+|x|)^a\,dx.
$$
Therefore we have proved \eqref{1.18}.

Let us assume now   that $|n-n_\infty|<c (1+|x|)^{-\delta}$. Then using  \eqref{MCEst} we conclude that
$$
\int |\nabla u-in_\infty^{1/2}\frac{x}{|x|}u|^2\,\frac{dx}{1+|x|}\leq c\int|f|^2(1+|x|)^a\,dx,
$$
which is \eqref{1.188}.
\newline
Only the uniqueness remains to be proved. From
$$
\int |\nabla u-in^{1/2}\frac{x}{|x|}u|^2\,\frac{dx}{1+|x|}\leq C
$$
we get that
\beq
{\rm liminf} \int_{|x|=r} \lp |\nabla u-in^{1/2}\frac{x}{|x|}u|^2\rp \  d\sigma(x) = 0, \qquad {\rm as} \quad r \to \infty,
\label{uniq1}
\eeq
Multiply in the equation $\Delta u+nu=0$ by $\bar u$ and integrate by parts in $|x|\leq R$. Taking the imaginary part one gets
${\cal I}m \int_{|x|=r}\bar u\partial_r u d\sigma(x) = 0.
$
And therefore
\beq
{\cal I}m \int_{|x|=R}\bar u(\partial_r u -in^{1/2}u )d\sigma(x) = -\int_{|x|=R}n^{1/2}| u|^2d\sigma(x).
\label{uniq2}
\eeq
Then from \eqref{MCEst}, \eqref{uniq1},  \eqref{uniq2}, and that $n$ is bounded we get that $u$ satisfies \eqref{Eq:0}. Therefore uniqueness follows from Theorem \ref{Thm0}.
\newline

\noindent {\bf Proof of Proposition \ref{Prop1}.} The argument is similar to the one to obtain \eqref{uniq1}. We multiply in
the equation $\Delta u+nu=f$ by $\bar u$, integrate by parts in the ball of radius $R_0$ and take
imaginary parts on both sides to obtain
\beq
{\cal I}m \int_{|x|=R_0}\bar ufd\sigma(x)={\cal I}m \int_{|x|=R_0}\bar u(\partial_r u -in^{1/2}u )d\sigma(x) = -\int_{|x|=R_0}n^{1/2}| u|^2d\sigma(x).
\label{prop3}
\eeq
 Integrate in $R_0$ in the above expression  for $0\leq R_0\leq R$, and divide by $R_0$. Recall that
  from \eqref{Eq:A7} $n\geq c_0\lambda$ for some $c_0>0$. Then the result follows from
Sommerfeld radiation condition and  \eqref{MCEst}. 

%
%%%%%%%%%%%%%%%%%%%%%%%%%%%%%%%%%%%%%%%%%%%%%%%%%%%
\section{Appendix}
\setcounter{equation}{0}
%%%%%%%%%%%%%%%%%%%%%%%%%%%%%%%%%%%%%%%%%%%%%%%%%%%
%%%%%%%%%%%%%%%%%%%%%%%%%%%%%%%%%%%%%%%%%%%%%%%%%%%
%
%
{\bf A.1.-  Basic identities.}
Our proof combines three basic
identities that have been used throughout this paper and that we
state here without proof (see \cite{PV} for a proof).
For real valued functions $\Psi, \; \vp,\; \psi \in {\cal S}(\R^d)$, we have
$$
    - \int_{\R^d}  \vp (x) | \nabla u(x) |^{2} + \frac{1}{2}
\int_{\R^d}  \Delta \vp(x) | u(x)|^{2}+ \int_{\R^d}  \vp(x) n(x) |
u(x)|^{2} \quad \quad 
$$
\beq
\quad \quad \quad \quad \quad \quad \quad \quad={\cal R}e  \int_{\R^d}  \vp(x)
\bu(x) f(x).
\label{Eq1}
\eeq
\medskip
$$
\e\int_{\R^d}  \psi(x) | u(x)|^{2}-{\cal I}m \int_{\R^d}  \nabla\psi (x)
\cdot \nabla u (x) \bu (x) \quad \quad  \quad \quad \quad \quad
$$
\beq
\quad \quad \quad \quad \quad \quad  \quad \quad = {\cal I}m \int_{\R^d}  f(x)
\bu (x)
\psi(x).
\label{Eq2}
\eeq
\medskip
$$
\int_{\R^d}  \lc\nabla \bu(x) \cdot D^{2}\Psi (x) \cdot \nabla u(x) -
\frac{1}{4} \Delta^{2}\Psi(x) | u(x)|^{2}+\frac{1}{2}
\nabla n(x)\cdot \nabla \Psi (x)| u(x)|^{2} \rc =
$$
\beq
-{\cal R}e \int_{\R^d}  f(x)\lp \nabla \Psi(x)\cdot \nabla \bu(x) +\frac{1}{2}
\Delta \Psi(x) \bu(x)\rp -
\e {\cal I}m\int_{\R^d}  \nabla \Psi (x) \cdot\nabla \bu (x) u(x).
\label{Eq3}
\eeq
\newline

{\bf A.2.- Some Examples.}
In this section we shall give examples of indices $n$ which satisfy
\beq\frac{(x\cdot \nabla  n(x))_{-}}{ n(x)}:=\tilde\beta < +\infty
\label{5.1}\eeq
with $\tilde\beta$ as small as wanted and such that Theorem
\ref{Thm1} does not hold true. In particular the estimate for the tangential 
derivative is false. Also we will exhibit the corresponding Sommerfeld radiation condition.

   Condition \fer{5.1} is weaker than \fer{Eq:A2} and appears
naturally in the study of the absence of embedded eigenvalues in the
continuous spectrum for the Schr\"odinger operator $\Delta +n$.
Recall at this respect the well known example due to Von Neumann and Wigner
of a potential which satisfies \fer{5.1} for $\tilde\beta$ large enough and
has an embedded eigenvalue, see \cite{RS} p. 233.

We give examples of wave guides which satisfy
\fer{5.1} but with a scaling which does not leave invariant that
condition, and therefore there is no possible
$\tilde\beta$ good for all of them. For these examples condition
\fer{Eq:A2} is not fulfilled either.

Define $Q$ as the unique positive solution with $Q(\pm\infty)=0$ of
$$
Q^{''}+\lp Q^2-\f{1}{2}\rp Q=0, \quad y\in \R.
$$
That is to say $Q(y)={sech} \ (y/\sqrt2)$. Also for $\lambda>0$ take
$Q_\lambda(y)=Q(\lambda y)$, which solves
$$
Q_\lambda^{''}+\lambda^2\lp Q_\lambda^2-\f{1}{2}\rp Q_\lambda=0.
$$
Set $\theta\in \cal C^{\infty}$ a bump function around the origin with
$\theta(x)=0$ if $|x|<1$ and $\theta(x)=1$ if $|x|>2$. Then call
$$u
^\e_\lambda(x,y)=Q_\lambda(y)\theta(x)e^{i\sqrt{1+i\e}|x|},
$$
with $0<\e<1$ and ${\cal I}m\sqrt{1+i\e}>0.$ Then $u^\e_\lambda$ solves
$$
\Delta u^\e_\lambda+\lp n_\lambda+i\e\rp u^\e_\lambda=f_\e\qquad (x,y)\in
{\R}^2,
$$
with
$$
n_\lambda:=\tilde n_\lambda+1-\f{1}{2}\lambda^2;\qquad \tilde
n_\lambda:=\lambda^2Q^2(\lambda y),
$$
and
$$f_\e(x,y)=\lp2i\sqrt{1+i\e} \ sig(x)
\theta^{'}(x)+\theta^{''}(x)\rp Q_\lambda(y)e^{i\sqrt{1+i\e}|x|}.
$$
Now it is straightforward to check that $N(f_\e)<\infty$. We can pass to the limit
in $\e$ and we get that, setting $u_\lambda=\lim_{\e\rightarrow 0^+}
u_\lambda^\e=Q_{\lambda}(y)\theta(x)e^{i|x|},$ then
$$
\Delta u_\lambda+n_{\lambda}u_\lambda=f(x,y):=\lp sig(x)
\theta^{'}(x)+\theta^{''}(x)\rp Q_\lambda(y)e^{i|x|}.
$$
Also if $0<\lambda<1/2$ and
$n_\infty^\lambda=1-\f{\lambda^2}{2}$,
we get
$$
n_\lambda-n_\infty^\lambda=\tilde n_\lambda=\lambda{^2}Q^2(\lambda y).
$$
Therefore
$$
r\f{\lp \dr\tilde  n_\lambda\rp_- }{n_\lambda}<2C\lambda^2\qquad r>1,
$$
and can be made as small as wanted.
On the other hand straightforward computations prove that $\frac{1}{R} \intbr \ |u_\lambda|^2<\infty$,
$$
\sup_\e\int |\nabla(e^{-i\sqrt{1+i\e}|x|}u^\e_\lambda|^2\leq c(1+\int(Q')^2\,dy)\le c,
$$
and that there is $c_0>0$ independent of $\e$ such that
$$
\int|\nabla_\tau u|^2\f{dx}{(1+(x^2+y^2)^{1/2}}\geq c_0\,| \log\e|.
$$
Also notice that $\varphi=i\sqrt{1+i\e}|x|$ is not a solution of $|\nabla\varphi|^2=n$.

%
%%%%%%%%%%%%%%%%%%%%%%%%%%%%%%%%%%%%%%%%%%%%%%%%%%%

\end{document}